\documentclass[11pt, a4paper]{article}

\usepackage[english]{babel}

\usepackage{amsmath}
\usepackage{amssymb}
\usepackage{amsthm}
\usepackage{tikz}
\usepackage{tikz-cd}
\usepackage{todonotes}
\usepackage{hyperref}
\usepackage{cleveref}
\usepackage{footmisc}
\usepackage{bbm}
\usepackage{afterpage}
\usepackage{mathrsfs}
\usepackage{sseq}
\usepackage{extarrows}
\usepackage{enumitem}
\usepackage{stmaryrd}
\usepackage{pdflscape}

\newcommand{\Cont}{\mathrm{Cont}}

\newcommand{\Ell}{\mathrm{Ell}}

\newcommand{\mult}{\mathrm{mult}}
\newcommand{\sm}{\mathrm{sm}}
\newcommand{\ord}{\mathrm{ord}}

\renewcommand{\top}{\mathrm{top}}
\newcommand{\Tate}{\mathrm{Tate}}

\newcommand{\BTtwo}{{{\mathrm{BT}^p_2}}}
\DeclareMathOperator{\Spec}{Spec}
\DeclareMathOperator{\Spf}{Spf}

\DeclareMathOperator{\KU}{K}
\DeclareMathOperator{\K}{K}
\renewcommand{\L}{\mathrm{L}}
\DeclareMathOperator{\KO}{KO}
\newcommand{\cKO}{{\KO^\wedge_\ast \hspace{-0.04cm}}}
\newcommand{\zKO}{{\KO^\wedge_0 \hspace{-0.04cm}}}

\newcommand{\ftKO}{{\KO^\wedge_{4t} \hspace{-0.02cm}}}
\newcommand{\cK}{{\K^\wedge_\ast \hspace{-0.04cm}}}
\newcommand{\zK}{{\K^\wedge_0 \hspace{-0.04cm}}}
\newcommand{\tK}{{\K^\wedge_t \hspace{-0.04cm}}}
\newcommand{\cE}{{E^\wedge_\ast \hspace{-0.04cm}}}

\newcommand{\cL}{{\L^\wedge_\ast \hspace{-0.04cm}}}
\newcommand{\zL}{{\L^\wedge_0 \hspace{-0.04cm}}}
\newcommand{\tsL}{{\L^\wedge_{t-s} \hspace{-0.02cm}}}

\newcommand{\smKO}{{\KO(\!(q)\!)}}
\newcommand{\smK}{{\K(\!(q)\!)}}

\newcommand{\Sph}{\mathbf{S}}

\DeclareMathOperator{\TMF}{TMF}
\DeclareMathOperator{\tmf}{tmf}
\DeclareMathOperator{\Tmf}{Tmf}

\DeclareMathOperator{\Alg}{Alg}
\DeclareMathOperator{\CAlg}{CAlg}

\DeclareMathOperator{\fDM}{fDM}

\DeclareMathOperator{\Mod}{Mod}

\newcommand{\id}{\mathrm{id}}

\newcommand{\CC}{\mathcal{C}}
\newcommand{\DD}{\mathcal{D}}
\newcommand{\E}{\mathbf{E}}

\newcommand{\F}{\mathbf{F}}
\newcommand{\G}{\mathbf{G}}
\newcommand{\h}{\mathrm{h}}

\newcommand{\M}{\mathcal{M}}
\newcommand{\N}{\mathbf{N}}

\renewcommand{\O}{\mathcal{O}}

\newcommand{\Q}{\mathbf{Q}}

\renewcommand{\P}{\mathbf{P}}

\newcommand{\T}{\mathrm{T}}

\newcommand{\x}{\mathfrak{X}}

\newcommand{\y}{\mathfrak{Y}}
\newcommand{\Z}{\mathbf{Z}}

\newcommand{\al}{\alpha}

\newcommand{\ga}{\gamma}
\newcommand{\Ga}{\Gamma}

\usepackage{mathtools}

\theoremstyle{theorem}\numberwithin{equation}{section}
\newtheorem{theorem}[equation]{Theorem}
\crefname{theorem}{{th}.\!\!}{{ths}.\!\!}
\Crefname{theorem}{{Th}.\!\!}{{Ths}.\!\!}
\newtheorem{theoremalph}{Theorem}

\crefname{theoremalph}{{th}.\!\!}{{ths}.\!\!}
\Crefname{theoremalph}{{Th}.\!\!}{{Ths}.\!\!}

\Crefname{problem}{{Prb}.\!\!}{{Prbs}.\!\!}
\newtheorem{prop}[equation]{Proposition}
\Crefname{prop}{{Pr}.\!\!}{{Prs}.\!\!}

\Crefname{lemma}{{Lm}.\!\!}{{Lms}.\!\!}
\newtheorem{cor}[equation]{Corollary}
\Crefname{cor}{{Cor}.\!\!}{{Cors}.\!\!}

\Crefname{conjecture}{{Conj}.\!\!}{{Conjs}.\!\!}

\theoremstyle{definition}\numberwithin{equation}{section}
\newtheorem{mydef}[equation]{Definition}
\Crefname{mydef}{{Df}.\!\!}{{Dfs}.\!\!}

\Crefname{recall}{{Rcl}.\!\!}{{Rcls}.\!\!}

\Crefname{construction}{{Con}.\!\!}{{Cons}.\!\!}

\Crefname{ass}{{As}.\!\!}{{As}.\!\!}

\Crefname{notation}{{Nt}.\!\!}{{Nts}.\!\!}

\Crefname{situation}{{St}.\!\!}{{Sts}.\!\!}

\theoremstyle{remark}\numberwithin{equation}{section}
\newtheorem{example}[equation]{Example}
\Crefname{example}{{Ex}.\!\!}{{Exs}.\!\!}
\newtheorem{nonexample}[equation]{Non-example}
\Crefname{nonexample}{{NonEx}.\!\!}{{NonEx}.\!\!}

\Crefname{claim}{{Clm}.\!\!}{{Clms}.\!\!}
\newtheorem{remark}[equation]{Remark}
\Crefname{remark}{{Rmk}.\!\!}{{Rmks}.\!\!}

\Crefname{idea}{{Id}.\!\!}{{Ids}.\!\!}

\Crefname{warn}{{Warn}.\!\!}{{Warns}.\!\!}

\Crefname{figure}{{Fig.}\!\!}{{Figs.}\!\!}
\Crefname{footnote}{{Fn.}\!\!}{{Fn.}\!\!}

\Crefname{part}{{\textsection}\!\!}{{\textsection}\!\!}
\Crefname{chapter}{{\textsection}\!\!}{{\textsection}\!\!}
\Crefname{section}{{\textsection}\!\!}{{\textsection}\!\!}
\Crefname{subsection}{{\textsection}\!\!}{{\textsection}\!\!}
\Crefname{appendix}{{\textsection}\!\!}{{\textsection}\!\!}

\begin{document}
\title{Uniqueness of real ring spectra up to higher homotopy}
\author{Jack Morgan Davies\footnote{\url{davies@math.uni-bonn.de}}}
\maketitle

\begin{abstract}
We discuss a notion of \emph{uniqueness up to $n$-homotopy} and study examples from stable homotopy theory. In particular, we show that the $q$-expansion map from elliptic cohomology to topological $K$-theory is unique up to $3$-homotopy, away from the prime $2$, and that upon taking $p$-completions and $\F_p^\times$-homotopy fixed points, this map is uniquely defined up to $(2p-3)$-homotopy. Using this, we prove new relationships between Adams operations on connective and dualisable topological modular forms---other applications, including a construction of a connective model of Behrens' $Q(N)$ spectra away from $2N$, will be explored elsewhere. The technical tool facilitating this uniqueness is a variant of the Goerss--Hopkins obstruction theory for \emph{real} spectra, which applies to various elliptic cohomology and topological $K$-theories with a trivial complex conjugation action as well as some of their homotopy fixed points.
\end{abstract}

\setcounter{tocdepth}{3}
\tableofcontents

\newpage

\addcontentsline{toc}{section}{Introduction}
\section*{Introduction}

Inside a $1$-category $\CC$, objects with a universal property are unique up to unique isomorphism, meaning the mapping set between any two such objects contains a single element. If $\CC=\h_1\DD$ is the homotopy $1$-category of an $\infty$-category $\DD$, then uniqueness in $\CC$ up to unique isomorphism translates to uniqueness in $\DD$ \emph{up to $1$-homotopy}. More generally, we say an object in $\DD$ with a property $P$ is \emph{unique up to $n$-homotopy} if all mapping spaces between any two objects with this property $P$ are $(n-1)$-connected, meaning their homotopy groups $\pi_d$ vanish in degrees $0\leq d\leq n-1$. In theory, studying the uniqueness of objects up to $n$-homotopy is most elegant for $n=\infty$, however, in practice, there are many objects of interest inside an $\infty$-category without an obvious or viable universal property.\\

In this article, we are interested in $\E_\infty$-rings of low chromatic height such as elliptic cohomology and topological $K$-theories, and our goal is to prove that certain morphisms between such $\E_\infty$-rings are unique up to $n$-homotopy for some finite $n$. To state our first theorem, let us write $\tmf$ for the $\E_\infty$-ring of \emph{connective topological modular forms} and $\KO\llbracket q\rrbracket$ for \emph{real Tate $K$-theory}; we recall these definitions in \Cref{backgoundandnotationsection}.

\begin{theoremalph}\label{mainuniquenessstatement}
A map of $\E_\infty$-rings $\tmf \to \KO\llbracket q\rrbracket[\frac{1}{2}]$ is uniquely determined by its effect upon applying complex $K$-theory homology in degree $0$, up to $3$-homotopy.\footnote{As one will see in our proof of \Cref{mainuniquenessstatement}, the $\tmf$ in \Cref{mainuniquenessstatement} can be replaced with $\TMF$ (but not with $\Tmf$ using our techniques; see \Cref{fromtmftotmf}), and the $\KO\llbracket q\rrbracket [\frac{1}{2}]$ with $\smKO[\frac{1}{2}]$ or $\KO[\frac{1}{2}]$; see \Cref{generalisedmainuniquenessstatmet} for more $p$-adic generalisations.}
\end{theoremalph}

For example, the topological $q$-expansion map $\tmf\to \KO\llbracket q\rrbracket$, itself a map of $\E_\infty$-rings, is uniquely determined up to $3$-homotopy, away from the prime $2$.\footnote{It is crucial that we work away from the prime $2$ in this article. At the end of the day, this is because the homotopy groups of $\KO_2$ are not supported in degrees divisible by $4$; see \Cref{failiureoftwopimraryremark}.} This theorem generalises the usual statement that such morphisms between the $\K(1)$-localisations of these $\E_\infty$-rings are unique up to $1$-homotopy; see \cite[Pr.A.6]{hilllawson} for example. We also show that when completed at an odd prime $p$, the above theorem admits generalisations to a large class of elliptic cohomology and topological $K$-theories, such as $\Tmf$, $\TMF_0(N)$, and $\KO$; see \Cref{generalisedmainuniquenessstatmet}. In another direction, if we complete at a prime $p$ and replace $\tmf_p$ and $\KO\llbracket q\rrbracket_p$ with their $G$-homotopy fixed points for a nontrivial finite subgroup $G\leq \F_p^\times$, where $\F_p^\times$ acts on these $\E_\infty$-rings through stable $p$-adic Adams operations, we obtain a stronger uniqueness statement depending on the order of $G$.

\begin{theoremalph}\label{mainuniquenessstatementpgrowth}
Let $p$ be an odd prime and $G\leq \Z_p^\times$ be a nontrivial finite subgroup with order $g$. A map of $\E_\infty$-rings $\tmf_p^{hG} \to \K\llbracket q\rrbracket_p^{hG}$ is uniquely determined by its effect on the zeroth $p$-adic complex $K$-theory homology group up to $(2g-1)$-homotopy.
\end{theoremalph}

When $G=\{\pm 1\}$, this is a $p$-complete version of \Cref{mainuniquenessstatement}. On the other extreme, if $G=\F_p^\times$ is the maximal finite subgroup, we obtain a uniqueness statement up to $(2p-3)$-homotopy. In this latter case, the $\E_\infty$-ring $\K\llbracket q\rrbracket_p^{hG}$ is a Tate $K$-theory analogue of the Adams summand $\L=\K^{h\F_p^\times}$ and $\tmf_p^{h\F_p^\times}$ is the \emph{height $2$ Adams summand} of \cite[\textsection3.1]{adamsontmf}.\\

Both \Cref{mainuniquenessstatement,mainuniquenessstatementpgrowth} are proven with applications in mind. Indeed, we treat these theorems almost as a ``$n$-homotopy extension property'' for the localisation map $\Tmf\to \TMF$, either away from $2$ or completed at an odd prime. To be precise, we use the above statements to prove that there is a family of homotopies between the stable Adams operations $\psi^k$ acting on $\tmf$ from \cite[Th.A]{adamsontmf}:
\[\psi^k\psi^\ell\simeq \psi^{k\ell}\colon \tmf[\frac{1}{k\ell}]\to \tmf[\frac{1}{k\ell}]\]
Moreover, this family of $1$-homotopies of morphisms of $\E_\infty$-rings can be chosen to be associative and unital up to $2$-homotopy. This is succinctly captured by the following result.

\begin{theoremalph}\label{adamsoperationshomotopiesintegral}
There is a homotopy of morphisms of $\E_\infty$-rings $\psi^{-1}\simeq \id\colon \tmf[\frac{1}{2}]\to \tmf[\frac{1}{2}]$ which is associative up to $2$-homotopy. More generally, fixing an integer $n$, there is a functor
\[\Psi^n\colon B(\prod_{p|n}\N)\to \h_2\CAlg\]
sending the unique point of $B(\prod_{p|n}\N)$ to $\tmf[\frac{1}{2n}]$ and $(p_1^{e_1},\ldots,p_r^{e_r})$ to $\psi^{\prod_i p_i^{e_i}}$, where $\h_2\CAlg$ denotes the homotopy $2$-category of $\E_\infty$-rings.
\end{theoremalph}

Using the techniques of \cite{adamsontmf}, where the operations $\psi^k$ are constructed, it is not clear that $\psi^{-1}\simeq \id$, let alone that the functor analogous to $\Psi^n$ into the homotopy $1$-category of $\E_\infty$-rings even exists. If we ``invert $\Delta$'' and replace $\tmf$ with the spectrum of periodic topological modular forms $\TMF$, then the corresponding statements for Adams operations can be found in \cite[Th.D]{heckeontmf}---in fact, these results for $\TMF$ are crucial in our proof of \Cref{adamsoperationshomotopiesintegral}. Other applications of \Cref{mainuniquenessstatement} include the construction of connective models of Behrens' $Q(N)$ spectra as well as operations thereupon; this is current work-in-progress. Such a construction was the original inspiration for this article.\\

Of course, there is also a $p$-adic version of \Cref{adamsoperationshomotopiesintegral} depending on a choice of $G\leq \Z_p^\times$.

\begin{theoremalph}\label{adamsoperationshomotopiespadic}
Fix an odd prime $p$ and $G\leq \Z_p^\times$ a nontrivial finite subgroup of order $g$. Then there is a functor
\[\Psi_{p,G}\colon B\Z_p^\times/G\to \h_{2g-2}\CAlg\]
sending the unique point of $B\Z_p^\times/G$ to $\tmf_p^{hG}$ and $\overline{k}\in \Z_p^\times/G$ to $\psi^k$.
\end{theoremalph}

Setting $G=\F_p^\times$ so that $g=p-1$, this theorem and \Cref{mainuniquenessstatementpgrowth} support the moral that stable homotopy theory, especially chromatic homotopy theory, becomes more sparse as the prime $p$ grows. \\

To prove all of the above theorems we use a variant of Goerss--Hopkins obstruction theory at odd primes $p$. Classically, one would use $p$-complete complex $K$-theory as the base cohomology theory to compute obstruction groups in the category of $\K(1)$-local $\E_\infty$-rings. If we did this we would only obtain uniqueness in \Cref{mainuniquenessstatement} up to $1$-homotopy, which is well-known, and the rest of our applications would falter. To amend this, we use a variant of Goerss--Hopkins obstruction theory now based on $\KU^{hG}$ for a nontrivial finite subgroup $G\leq \Z_p^\times$ following an idea of Behrens at the prime $p=2$. The following is the analogue of the results of the appendix \cite[\textsection 12.A]{tmfbook}; see \Cref{bottperiodiceinftyringdefinition} for the definition of a \emph{real} $\E_\infty$-ring.

\begin{theoremalph}\label{ghtheorymain}
Fix an odd prime $p$, a nontrivial finite subgroup $G\leq \Z_p^\times$, and two $p$-complete real $\E_\infty$-rings $A,B$ where $B$ is $\K(1)$-local. Let us also write $\KU^{hG}=\L$. Then for any map of $\E_\infty$-rings $f\colon A\to B$ there exists a second quadrant spectral sequence converging to $\pi_{t-s} (\CAlg(A, B), f)$ with $E_2$-page given as
\[E_2^{0,0}=\Alg^{\bar{\theta}}_{\L_\ast}(\cL A,\cL B)\]
\[E_2^{s,t}=H^s_{\bar{\theta}}(\cL A/\L_\ast,\Omega^t \cL B)\qquad t\geq 1\]
and zero elsewhere.
\end{theoremalph}

The fact that both $\TMF$ and $\KO\llbracket q\rrbracket$ are real, essentially because the Adams operation $\psi^{-1}$ acts trivially on these spectra, allows us to use the above theorem combined with a vanishing statement to prove a $p$-adic version of \Cref{mainuniquenessstatement}. More generally, we use an algebro-geometric criterion to show that many of our favourite elliptic cohomology theories and topological $K$-theories are real; see \Cref{hypothesesarefulfilledone}.

\subsection*{Outline}

This article has three sections: background, construction of our obstruction theory with real $\theta$-algebras and real spectra, and the applications of this obstruction theory to \Cref{mainuniquenessstatement} and \Cref{adamsoperationshomotopiesintegral}. In more detail:

\begin{itemize}
\item In \Cref{uniquenesssection}, the concept of \emph{uniqueness up to $n$-homotopy}, in the sense we will use it, is formally defined and discussed.
\item In \Cref{backgoundandnotationsection}, we recall definitions from $\K(1)$-local stable homotopy theory such as $\K(1)$-local homology theories, $\theta$-algebras, and their reduced variants.
\item In \Cref{comodulesforko}, we define and discuss real $\theta$-algebras and real reduced $\theta$-algebras, which we call real $\bar{\theta}$-algebras. In particular, we show that the categories of graded real $\theta$-algebras and graded real $\bar{\theta}$-algebras are equivalent and compare their Andr\'{e}--Quillen cohomology.
\item In \Cref{kohomoologyofellcohom}, the notion of \emph{real} spectra is introduced, which allows us to prove \Cref{ghtheorymain}.
\item In \Cref{criterionforbpsection}, we discuss an algebro-geometric condition which implies that the certain elliptic cohomology or topological $K$-theories are real, meaning $\{\pm1\}$-real. These techniques are then also used to show various other $G$-homotopy fixed point spectra are real.
\item In \Cref{uniquenessofqexpasniosection}, we prove a vanishing statement for the Andr\'{e}--Quillen cohomology of $\bar{\theta}$-algebras. A generalisation of \Cref{mainuniquenessstatement} is then proven as a corollary of this vanishing statement together with \Cref{ghtheorymain}, along with \Cref{mainuniquenessstatementpgrowth}.
\item In \Cref{adamsoperationssection}, the $\E_\infty$-ring $\Tmf$ is decomposed into the periodic spectra $\TMF$ and $\KO\llbracket q\rrbracket$ glued together along $\smKO$. The $1$- and $2$-homotopies between Adams operations on $\TMF$ and $\KO\llbracket q\rrbracket$ (using their modular description) are then glued together on $\smKO$ using \Cref{mainuniquenessstatementpgrowth} which yields \Cref{adamsoperationshomotopiespadic}. We then obtain \Cref{adamsoperationshomotopiesintegral} by further gluing in some rational information.
\end{itemize}

\subsection*{Notation}

Throughout this article, $p$ will be an odd\footnote{Many of the statements of this article hold for $p=2$, however, the vanishing statement \Cref{vanishingofobstructions} at $p=2$ has to be altered to reflect the $2$-local homotopy groups of $\KO$; see \cite[Lm.12.7.13]{tmfbook}. This alteration then ruins the prime $2$ versions of all of our results from the introduction. For this reason, as well as the facts that many arguments are streamlined at odd primes and Behrens has made many of our statements here at $p=2$ (see \cite[\textsection 12.7 \& \textsection 12.A]{tmfbook}), we restrict this whole article to odd primes $p$.} prime, $G$ will denote a nontrivial finite subgroup of $\Z_p^\times$ which we often think of as a subgroup of the groups of lifts $\F_p^\times\leq \Z_p^\times$, and $\Ga=\Z_p^\times/G$ will denote the quotient group. We will also use $g$ to denote the order of $G$. Let $\KU$ be the $\E_\infty$-ring of \emph{($p$-complete) periodic complex $K$-theory}, and write $\L=\KU^{hG}$ for the $G$-homotopy fixed points of $\KU$ with respect to the stable $p$-adic Adams operations. In particular, when $G=\{\pm 1\}$ we have $\L=\KO$, and when $G=\F_p^\times$ this is the classical periodic Adams summand. Both $\KU$ and $\L$ have $\E_\infty$-actions of $\Z_p^\times$ and $\Z_p^\times/G=\Ga$, respectively, given by these Adams operations $\psi^k$; see \cite[\textsection5.5]{luriestheorem} for example. 

\subsection*{Acknowledgements}

Thank you to Venkata Sai Narayana Bavisetty for stimulating discussions about $Q(N)$ which eventually led to this project. I’d like to also thank Tommy Lundemo and Lennart Meier for their stimulating conversations and feedback.


\section{Definitions and notation}

We will freely use the language of $\infty$-categories to ground our homotopy theory.


\subsection{Uniqueness up to $n$-homotopy}\label{uniquenesssection}

The study of uniqueness in homotopy theory here will be quite superficial and specialised to keep our discussion here elementary.

\begin{mydef}
Let $f\colon X\to Y$ be a map of spaces and $y$ be a point in $Y$. Write $M_y=X\times_Y \{y\}$ for the \emph{moduli space of lifts of $y$}. For a lift $x$ of $y$, so an element $x\in M_y$, we say $x$ is \emph{unique up to $1$-homotopy} if $M_y$ is connected, ie, $\pi_0 M_y$ is the singleton set $\{x\}$. For $2\leq n\leq \infty$, we say that a unique lift $x$ of $y$ up to $1$-homotopy is \emph{unique up to $n$-homotopy} if the group $\pi_k(M_y,x)$ vanishes for all $0<k<n$.\footnote{If $x$ is unique up to $\infty$-homotopy, one often says that $x$ is \emph{unique up to contractible choice}.}
\end{mydef}

The above definition is motivated by a simple observation: say we have a unique lift $x$ of $y$ up to $1$-homotopy, and let $x'$ be another lift of $y$. This uniqueness of $x$ yields a homotopy $h\colon \Delta^1\to M_y$ from $x$ to $x'$. To quantify if this homotopy $h$ recognising the uniqueness of $x$ is unique up to $2$-homotopy, take another homotopy $h'\colon \Delta^1\to M_y$ from $x$ to $x'$. The composite path $h'\star h^{-1}\colon \Delta^1\to M_y$ determines a loop $\ga$ in $M_y$ based at $x$. If $x$ is unique up to $2$-homotopy, then $\pi_1(M_y,x)$ vanishes and this loop $\ga$ can be contracted to the constant loop at $x$ by a $2$-homotopy $H\colon \Delta^2\to M_y$. This $2$-homotopy $H$ also takes the form of a $2$-homotopy between $h$ and $h'$, witnessing that homotopies recognising the uniqueness of $x$ are homotopic. The story goes on, with uniqueness up to $n$-homotopy meaning that potentially different $(n-1)$-homotopies can be glued together to form a map $\partial \Delta^{n+1}\to M_y$, which can be lifted to a map from $\Delta^{n+1}$ if $\pi_n(M_y,x)$ vanishes.\\

The reader is encouraged to consider examples such as the uniqueness of maps of pointed space $S^n\to S^n$ recognising multiplication by $2$ on $\pi_n$ for various $n$. For a more complicated example also relevant to this article, it is shown in \cite[Th.A]{uniqueotop} that the sheaf $\O^\top$ defining $\Tmf$ is uniquely determined up to $1$-homotopy by the fact that it defines natural elliptic cohomology theories on the small affine étale site of ${\M}_\Ell$.\footnote{We believe that $\O^\top[\frac{1}{2}]$ is uniquely determined by this property up to $3$-homotopy by combining the ideas of \cite{uniqueotop} together with \Cref{mainuniquenessstatement}.}\\

The following is a simple criterion for calculating homotopical uniqueness.

\begin{prop}\label{reductiontohomotopygroups}
    Let $f\colon X\to Y$ be a map of spaces and $y$ be an element of $Y$ such that the connected component of $y$ is contractible, for example, if $Y$ itself is discrete. Then a unique lift $x$ of $y$ up to $1$-homotopy is unique up to $n$-homotopy if and only if $\pi_k(X,x)$ vanishes for $0<k<n$.
\end{prop}

\begin{proof}
Writing $Y_y$ for the connected component of $y$ and $X_y=X\times_Y Y_y$, notice there is a natural equivalence of spaces $M_f\simeq X_y\times_{Y_y}\{y\}$ from the pasting lemma for Cartesian diagrams. This leads to the fibration $M_y\to X_y \to Y_y$ over the point $y$, now viewed as an element in $Y_y$. The desired result now follows from the long exact sequence of the above fibration.    
\end{proof}


\subsection{$\K(1)$-local homology theories and comodules}\label{backgoundandnotationsection}

Fix an odd prime $p$ and implicitly work $\K(1)$-locally. In particular, our algebra will be completed at a prime $p$, hence tensor products will be $p$-completed over $\Z_p$. We will also write stacks with a subscript $\Z_p$ which indicates the Cartesian product with the formal stack $\Spf\Z_p$ (as opposed to $\Spec \Z_p$).

\begin{mydef}
Let $E$ and $X$ be spectra. Define the \emph{($\K(1)$-local) $E$-homology of $X$} by the formula
\[\cE X=\pi_\ast L(E\otimes X)\]
where $L$ is localisation at the first Morava $K$-theory $\K(1)$.
\end{mydef}

This is not a homology theory in the classical sense, as infinite direct sums are not preserved, however, it is the sensible replacement when working with the $\K(1)$-local stable homotopy category; see \cite[\textsection9]{hoveysticklandMoravalocal}. For spectra $X$ such that $\cK X$ is torsion-free, it follows that $\cK X$, which \emph{a priori} is only $L$-complete, is actually $p$-complete. As we will often make this torsion-free assumption, we will sometimes call $\cK X$ the $p$-adic $\KU$-homology of $X$.\\

The natural home for our $\KU$- and $\L$-homology groups are the following algebraic categories.

\begin{mydef}
A \emph{graded Morava module} or \emph{graded $\psi$-module} is a $\Z$-graded $\K_\ast$-module $M_\ast$ such that each $M_n$ is $L$-complete together with a continuous $\Z_p^\times$-action, which we denote as Adams operations $\psi^k$ for $k\in\Z_p^\times$, such that for all $a\in \K_\ast$ and $m\in M_\ast$ we have $\psi^k(am)=\psi^k(a)\psi^k(m)$ for all $k\in \Z_p^\times$. There is a related ungraded notion over $\Z_p$, where the action of $\Z_p^\times$ on $\Z_p$ is trivial. Denote the categories of $\psi$-modules (resp.\ graded $\psi$-modules) by $\Mod^{\psi}_{\Z_p}$ (resp.\ $\Mod^{\psi}_{\K_\ast}$). For a fixed nontrivial finite subgroup $G\leq \Z_p^\times$, a \emph{reduced graded Morava module} or \emph{graded $\bar{\psi}$-module} is a $\Z$-graded $L$-complete $\L_\ast$-module $M_\ast$ with a continuous $\Ga$-action with the same compatibility as $\psi$-modules, with $\L$ replacing $\K$. There is also an obvious ungraded notion of $\bar{\psi}$-module over $\Z_p$. Write $\Mod_{\Z_p}^{\bar{\psi}}$ and $\Mod_{\L_\ast}^{\bar{\psi}}$ for these categories of $\bar{\psi}$-modules.
\end{mydef}

The appropriate multiplicative objects in these categories are $\theta$-algebras.

\begin{mydef}\label{thetaalgebradefinition}
A \emph{$\theta$-algebra} is a $\Z_p$-algebra $A_0$ equipped with the structure of a $\psi$-module and a $\Z_p^\times$-equivariant homogeneous operator $\theta$ with $\theta(1)=0$ such that the following two conditions hold for all $x,y\in A_0$:
\begin{enumerate}
\item $\theta(xy)=x^p\theta(y)+\theta(x)y^p+p\theta(x)\theta(y)$
\item $\theta(x+y)=\theta(x)+\theta(y)-\frac{1}{p}\sum_{i=1}^{p-1} \binom{p}{i}x^i y^{p-i}$
\end{enumerate}
There is also a graded version, whose definition we safely leave to the reader, as all graded $\theta$-algebras we will see in this article are of the form $\KU_\ast\otimes A_0$. Write $\Alg^{\theta}_{\K_\ast}$ and $\Alg^\theta_{\Z_p}$ for the categories of graded and ungraded $\theta$-algebras. A \emph{reduced graded $\theta$-algebra} is a $\Z$-graded $\L_\ast$-algebra $A_\ast$ equipped with the structure of a $\bar{\psi}$-module and a $\Ga$-equivariant operation $\theta$ satisfying the same axioms as above. Also define ungraded $\bar{\theta}$-algebras over $\Z_p$, and write $\Alg^{\bar{\theta}}_{\L_\ast}$ and $\Alg^{\bar{\theta}}_{\Z_p}$ for the categories of graded and ungraded $\bar{\theta}$-algebras.
\end{mydef}

The following crucial theorem connects these algebraic notions to homotopy theory. 

\begin{theorem}\label{wehavethetaalgerbsandsuch}
Let $X$ be a spectrum and $A$ be an $\E_\infty$-ring. Then $\cK X$ is a $\psi$-module, $\cK A$ is a $\theta$-algebra, $\cL X$ is a $\bar{\psi}$-module, and $\cL A$ is a $\bar{\theta}$-algebra.
\end{theorem}

\begin{proof}
We refer the reader to \cite[Pr.2.2.2 \& Th.2.2.4]{gh05} for the first two statements, and the second two follow as $\cL Y$ can be computed for any spectrum $Y$ as the homotopy groups of $(\K\otimes Y)^{hG}$ by $\K(1)$-local ambidexterity and a collapsing $G$-homotopy fixed point spectral sequence as $|G|||\F_p^\times|=p-1$ is invertible. Indeed, as discussed in \cite[Lm.12.7.8]{tmfbook} (where a reference to \cite[Th.1.5]{kuhnambi} is given), ambidexterity states that the natural map
\[Z_{hG} \xrightarrow{\simeq} Z^{hG}\]
is a $\K(1)$-equivalence for all spectra $Z$ with $G$-action, leading to the chain of equivalences
\[\L {\otimes} Y\xleftarrow{\simeq} \KU_{hG} {\otimes} Y\xleftarrow{\simeq} (\KU {\otimes}Y)_{hG}\xrightarrow{\simeq} (\KU {\otimes} Y)^{hG}\]
using the fact that the smash product is cocontinuous in each variable.
\end{proof}

The forgetful functor from $\theta$-algebras to $\psi$-modules has a left adjoint written as $\P_\theta$. We will also write $\P_{\theta,\ast}$ for the graded analogue, as well as $\P_{\bar{\theta}}$ and $\P_{\bar{\theta}, \ast}$ in the reduced cases. One can now use such ``free functors'' to construct simplicial resolutions and define Andr\'{e}--Quillen cohomology groups; see \cite[\textsection4]{gh04} or \cite[\textsection2.4]{gh05}. For us, given a graded $\theta$-algebra $A_\ast$ and a $\psi$-module $M_\ast$, we define $H^s_\theta(A_\ast/\K_\ast,M_\ast)$ as the cohomology of the cochain complex associated with the cosimplicial object
\[\Alg^{\theta}_{\K_\ast/A_\ast}(P_\bullet,A_\ast\ltimes M_\ast)\]
where $P_\bullet\to A_\ast$ is a simplicial resolution of $A_\ast$ by free graded $\theta$-algebras on $\psi$-modules which are projective as $\K_\ast$-modules. The same goes for $\bar{\theta}$-algebras.


\section{Real spectra and their Goerss--Hopkins obstruction theory}

The notion of a \emph{real spectrum} will be useful for us, as this strong condition will simultaneously allow us to reduce many statements about $\L$-homology to $\KU$-homology. When $G=\{\pm 1\}\leq \Z_p^\times$, then we will see in \Cref{criterionforbpsection} that many of our favourite elliptic cohomology and topological $K$-theories, such as $\KO$, are real with respect to $G=\{\pm 1\}$.

Fix an odd prime $p$, a nontrivial finite subgroup $G\leq \Z_p^\times$, and implicitly localise everything at $\K(1)$.


\subsection{Real $\psi$-modules and $\bar{\psi}$-modules}\label{comodulesforko}

The following notions of \emph{real} $\theta$-algebras and spectra are blatantly adapted from Behrens' definition \cite[Df.12.7.10]{tmfbook}.

\begin{mydef}
Given a graded $\psi$-module $M_\ast$, we say that $M_\ast$ is \emph{real} if the natural map of $\psi$-modules $\KU_\ast\otimes M_0\to M_\ast$ is an isomorphism and if the $G$-action on $M_0$ is trivial. Given a graded $\bar{\psi}$-module $N_\ast$, we say that $N_\ast$ is \emph{real} if the natural map $\L_\ast\otimes N_0\to N_\ast$ is an equivalence. If we have to clarify which $G$ the adjective real refers to, we will write \emph{real with respect to $G$}.
\end{mydef}

The justification for this definition is given by the following theorem. First, note that there is a functor
\[F\colon \Mod^{\bar{\psi}}_{\L_\ast}\to \Mod^{\psi}_{\KU_\ast}\qquad M_\ast\mapsto \KU_\ast\otimes_{\L_\ast} M_\ast\]
where the $\Z_p^\times$-action on $F M_\ast$ can be described by noting that $\L\to \KU$ is $\Z_p^\times$-equivariant, where the action on $\L$ is restricted from its natural $\Ga$-action. 

\begin{theorem}\label{equivalencesofrealness}
Fix an odd prime $p$ and a nontrivial finite subgroup $G\leq \Z_p^\times$. The functors
\[F\colon \Mod^{\bar{\psi}}_{\L_\ast}\to \Mod^{\psi}_{\KU_\ast}\qquad M_\ast\mapsto \KU_\ast\otimes_{\L_\ast} M_\ast=L M_\ast\]
\[\{\pm 1\}\colon\Mod^{\psi}_{\KU_\ast}\to \Mod^{\bar{\psi}}_{\L_\ast}\qquad N_\ast\mapsto N_\ast^{\{\pm1\}}\]
restricts to an equivalence between real $\bar{\psi}$-modules and real $\psi$-modules. Moreover, both of these functors lift to $\Alg^{\bar{\theta}}_{\L_\ast}$ and $\Alg^{\theta}_{\L_\ast}$ through the evident forgetful functors, and hence also restrict to equivalences between real $\bar{\theta}$-algebras and real $\theta$-algebras.
\end{theorem}

\begin{proof}
The natural isomorphisms
\[(\K_\ast\otimes_{\L_\ast} M_\ast)^{G}\simeq (\K_\ast\otimes M_0)^{G}\simeq \L_\ast\otimes M_0\simeq M_\ast\]
\[\K_\ast\otimes_{\L_\ast} N_\ast^{G}\simeq \K_\ast\otimes_{\L_\ast}(\K_\ast\otimes N_0)^{G}\simeq \K_\ast\otimes_{\L_\ast}(\L_\ast\otimes N_0)\simeq N_\ast\]
show these functors are mutual inverses when restricted to real modules. Above we have used that $G$-fixed points commute with $p$-completed tensor products, which follows as the additive norm $X_{G}\to X^{G}$ from coinvariants to invariants is an isomorphism as the order of $G$ is invertible from the divisibility relation $|G|||\F_p^\times|=p-1$.
\end{proof}

A key observation that we will use many times later, is that given a real spectrum $X$, then either $\cL X$ or $\cK X$ determine each other---an immediate consequence of \Cref{bottperiodicspectragivegood} and the above theorem. Another corollary of the above is that we can calculate Andr\'{e}--Quillen cohomology of a real $\bar{\theta}$-algebra $A_\ast$ as the Andr\'{e}--Quillen cohomology of $LA_\ast$ as a ${\theta}$-algebra.

\begin{cor}\label{cohomologyiseasyforrealguys}
Given a graded real $\bar{\theta}$-algebra $A_\ast$ and a graded real $\bar{\psi}$-module $M_\ast$, then the functor $F$ induces an isomorphism of abelian groups
\[H^s_{\bar{\theta}}(A_\ast/\L_\ast,M_\ast)\simeq H^s_{{\theta}}(FA_\ast/\K_\ast,FM_\ast),\qquad s\geq 0.\]
\end{cor}

\begin{proof}
The assumption that $A_\ast$ is real yields $A_\ast=\L_\ast\otimes A_0$, so let us choose a free resolution $P_\bullet$ of $A_\ast$ as a graded $\bar{\theta}$-algebra of the form $P_\bullet=\L_\ast\otimes P'_\bullet$ where $P'_\bullet\to A_0$ is a free simplicial $\bar{\theta}$-algebra resolution. Let us choose $P'_\bullet$ a little more carefully. As explained by Behrens' for the prime $p=2$ in \cite[\textsection 12.A]{tmfbook} (on the page above formula (A.2)), we can choose topological generators $\{x_\al\}$ of $A_0$ as a $\bar{\theta}$-algebra, and we may further take each $x_\al$ to have open isotropy inside $\Ga=\Z_p^\times/G$, else we choose other generators. Hence, for each $\al$ there is a $j\geq 1$ such that the generator $x_\al$ defines a map of $\bar{\theta}$-algebras
\[x_\al\colon \P_{\bar{\theta}}\left(\Z_p[(\Z/p^j\Z)^\times/G]\right)\to A_0\]
such that map of $\bar{\theta}$-algebras
\[P'_0=\P_{\bar{\theta}}\left({\bigoplus_{\al}} \Z_p[(\Z/p^j\Z)^\times/G]\right)\to A_0\]
defined by all of the generators, where the direct sum above is implicitly $L$-completed, is surjective. By restricting the $\Ga$-action on the $\bar{\theta}$-algebra above to a $\Z_p^\times$-action along the quotient of groups $\Z_p^\times\to \Ga$, we see that $P'_0$ is the start of a free simplicial resolution for $A_0$ as a $\theta$-algebra now. Indeed, both free functors $\P_{\theta}$ or $\P_{\bar{\theta}}$ simply adjoin a free operator $\theta$ commuting with the $\Z_p^\times$- or $\Ga$-action and satisfying the axioms of \Cref{thetaalgebradefinition}. These operators agree if $G$ acts trivially as it does on $P'_0$. Similarly, analysing the kernel of $P'_0\to A_0$ leads to $P'_1\to P'_0$, and inductively we obtain the whole free simplicial resolution $P'_\bullet$ of $A_0$, importantly, as \textbf{both} a $\theta$-algebra and a $\bar{\theta}$-algebra. In particular, as $F$ does not change alter degree zero, we see that $FP_\bullet\to FA_\ast$ is a free simplicial resolution of $A_\ast$ as a graded $\theta$-algebra. This leads to the isomorphism of cosimplicial objects
\[\Alg^{\bar{\theta}}_{\L_\ast/A_\ast} (P_\bullet,A_\ast\ltimes M_\ast) \xrightarrow{\simeq} \Alg^{\theta}_{\KU_\ast/F A_\ast}(F P_\bullet,FA_\ast\ltimes FM_\ast)\]
which upon taking associated cochain complexes and cohomology yields the desired isomorphism.
\end{proof}


\subsection{Real spectra and their Goerss--Hopkins spectral sequence (\Cref{ghtheorymain})}\label{kohomoologyofellcohom}

Our interest in real $\psi$-modules and $\bar{\psi}$-modules comes from our interest in the class of \emph{real} spectra.

\begin{mydef}\label{bottperiodiceinftyringdefinition}
For a given nontrivial finite subgroup $G\leq \Z_p^\times$, we say a spectrum $X$ is \emph{real} if $\cK X$ is torsion-free, concentrated in even degrees, and that the induced $G$-action in degree zero is trivial.
\end{mydef}

The case when $G=\{\pm 1\}$ is our motivation for this notation---a real spectrum $X$ with respect to $\{\pm1\}$ is one such that inclusion of fixed points $\zKO X\to \zK X$ given by the complex conjugation action is an isomorphism, as we will shortly see.\\

This definition only depends on the $\K(1)$-localisation of $X$ for a fixed (odd) prime $p$---we will not need a potential integral definition in this article. A real spectrum should be thought of as ``a nice spectrum with trivial $\psi^\ga$-action for $\ga\in G$''. For example, using \Cref{exampleofinversefixed} and \Cref{hypothesesarefulfilledone} we will see that $\tmf$, $\TMF_0(N)$, $\KO$, and $\KO\llbracket q\rrbracket$ are all real. Some non-examples include $\K$, as the classical calculation
\[\zK\K\simeq \Cont(\Z_p^\times,\Z_p)\]
shows that the $\{\pm1\}$-action on the left, translated to the conjugation action on the right, is nontrivial. The eager reader can similarly show that $\TMF_1(N)$ is \textbf{not} real for $N\geq 3$ and $G=\{\pm1\}$ using \Cref{ktheorycalculation} and showing the $\{\pm 1\}$-action there is nontrivial. Moreover, for $G\neq \{\pm 1\}$, one can also check that $\Tmf$ itself is not real, again using \Cref{ktheorycalculation}.\\

It follows rather easily that the $\L$-homology of a real $\E_\infty$-ring is real.

\begin{prop}\label{bottperiodicspectragivegood}
Let $X$ be a spectrum such that $\cK X$ is torsion-free and concentrated in even degrees. Then $X$ is real if and only if the map
\begin{equation}\label{complexifcationisomorphism}\zL X\to \zK X\end{equation}
induced by the inclusion of fixed points $\L\to \K$, is an isomorphism. If $X$ is real spectrum, then the natural maps\begin{equation}\label{kohomologycalcuation}\KU_\ast\otimes \zK X\xrightarrow{\simeq} \cK X\qquad\qquad \L_\ast\otimes \zL X\xrightarrow{\simeq} \cL X\end{equation}
are isomorphisms. In particular, the $\psi$-modules $\cK X$ and the $\bar{\psi}$-module $\cL X$ are both real in their respective categories.
\end{prop}

The condition that (\ref{complexifcationisomorphism}) was an isomorphism was Behrens' original definition of a real (what he called \emph{Bott periodic} for $G=\{\pm 1\}$ and at $p=2$), so the above can further be seen as a reconciliation between our two definitions.

\begin{proof}
Let us first note that the fact that $\cK X$ is torsion-free and concentrated in even degrees implies that the natural map
\begin{equation}\label{basechangeofcomplexktheory}\K_\ast\otimes\zK X\xrightarrow{\simeq} \cK X\end{equation}
is an isomorphism of graded $\theta$-modules. As another remark, notice that for any spectrum $Y$, then the natural map $\L\otimes Y\to \K\otimes Y$ is the inclusion of $\{\pm1\}$-homotopy fixed points by $\K(1)$-local ambidexterity; we also used this fact in the proof of \Cref{wehavethetaalgerbsandsuch}. In particular, we can now consider the following $G$-homotopy fixed point spectral sequence (HFPSS):
\begin{equation}\label{hfpssforx}E_2^{s,t}\simeq H^s(G,\tK X)\simeq H^s(G,\K_\ast\otimes \zK X)\Longrightarrow \tsL X\end{equation}
The above spectral sequence collapses immediately, as the groups $\cK X$ are $p$-complete for an odd prime $p$, hence the order of $G$ is invertible, and the edge map
\begin{equation}\label{edgemap}\cL X \xrightarrow{\simeq} (\cK X)^{G}\end{equation}
is an isomorphism. In degree zero the above map is the desired inclusion of fixed points map (\ref{complexifcationisomorphism}), and it is now clear that this map is an isomorphism if and only if the $G$-action on $\zK X$ is trivial.\\

The calculation of the $E_\infty$-page of the HFPSS (\ref{hfpssforx}) can be summarised by stating that there is a natural isomorphism of graded reduced $\theta$-modules
\[\L_\ast\otimes \zK X \xrightarrow{\simeq} \cL X\]
which when combined with (\ref{complexifcationisomorphism}) yields the second isomorphism of (\ref{kohomologycalcuation}); the first is (\ref{basechangeofcomplexktheory}).
\end{proof}

Essentially as a corollary of our study of ``reality'' above, we can prove \Cref{ghtheorymain}

\begin{proof}[Proof of \Cref{ghtheorymain}]
By \cite[Th.2.4.14]{gh05}, our hypotheses produce a spectral sequence converging to the homotopy groups of the space $\CAlg(A,B)$ based at $f$, so we are left to compute the $E_2$-page. Using \Cref{bottperiodicspectragivegood}, we see that $\cK X$ and $\cL X$ are both real if $X$ is real, so this applies to both $X=A,B$. For $s=t=0$, the given $E_2$-page takes the form
\[E_2^{0,0}=\Alg^\theta_{\KU_\ast}(\cK A,\cK B)\simeq \Alg^{\bar{\theta}}_{\L_\ast}(\cL A,\cL B)\]
using \Cref{equivalencesofrealness}, and for $t\geq 1$
\[E_2^{s,t}=H^s_\theta(\cK A/\K_\ast,\cK B)\simeq H^s_{\bar{\theta}}(\cL A/\L_\ast,\cL B)\]
using \Cref{cohomologyiseasyforrealguys}, as desired.
\end{proof}


\section{Uniqueness of maps of between real $\E_\infty$-rings}

In this section, we apply the obstruction theory of \Cref{ghtheorymain} to morphisms of certain $\E_\infty$-rings.


\subsection{Real elliptic cohomology and $K$-theories}\label{criterionforbpsection}

The condition that the $G$-action on $\zK X$ be trivial can naturally be tricky to determine, so it can be difficult to show that various spectra are real. If $X$ is an elliptic cohomology theory or a form of $K$-theory and $G=\{\pm 1\}$, then we have some hope. Write $\M_\Ell$ for the \emph{compactified moduli stack of elliptic curves}, so the stack which classifies generalised elliptic curves; see \cite{dr}.

\begin{mydef}
Let $f\colon \x\to \M_\Ell$ be a morphism of formal Deligne--Mumford stacks determined by a generalised elliptic curve $C$ over $\x$. Consider $f$ as an object in the slice $2$-category $\fDM_{/{\M}_\Ell}$. There is an involution $\tau$ of $f$ in this $2$-category defined by the pair $(\id_\x,[-1])$ where $[-1]\colon C\to C$ is the inversion isomorphism on the generalised elliptic curve $C$. We say that $f$ is \emph{inverse fixed} if there exists a trivialisation of $\tau$, so a $2$-morphism in $\fDM_{/{\M}_\Ell}$ between $\tau$ and the identity.
\end{mydef}

The following remarks can look tautological, to begin with, so the reader is advised to read the following example and non-example in parallel.

\begin{example}\label{exampleofinversefixed}
Consider $\x=\M_{\Ell,\Z_p}=\M_\Ell\times \Spf \Z_p$. In this case, there is a $2$-morphism in $\fDM_{/\M_\Ell}$ given by $[-1]\colon (\id,\id)\to (\id,[-1])=\tau$. In the case of $\x=\M^\sm_0(n)_{\Z_p}$ things are a little less tautological. Here, $\M^\sm_0(n)$ is the moduli stack of smooth elliptic curves equipped with a cyclic subgroup of order $n$ for a positive integer $n$ not divisible by $p$; see \cite[\textsection5]{km}. In this case, we can also define a $2$-morphism $[-1]\colon (\id,\id)\to(\id,[-1])$ as the $[-1]$-action of an elliptic curve fixes a choice of cyclic subgroup. The sections of the sheaf $\O^\top$ on these stacks are $\Tmf_p$ and $\TMF_0(N)_p$, respectively. The same $2$-morphism also works for $\x=\M_{\Tate,\Z_p}$, where $\M_\Tate$ is the Tate moduli stack, and also for the smooth Tate moduli stack $\M_\Tate^\sm$, and the moduli stack $\M_{\G_m}$ of forms of $\G_m$; see \cite[\textsection3.5]{hilllawson} or \cite[Dfs.1.1-2]{adamsontmf}. The $\E_\infty$-rings assocated to these stacks are $\KO\llbracket q\rrbracket_p$, $\smKO_p$, and $\KO_p$, respectively.
\end{example}

\begin{nonexample}
The moduli stack $\x=\M^\sm_1(n)_{\Z_p}$ is \textbf{not} inverse fixed for $n\geq 3$. Here $\M^\sm_1(n)$ is the moduli stack of smooth elliptic curves with a chosen point of order $n$ for some positive integer $n$ not divisible by $p$. To see this is not inverse fixed, one can use the converse to \Cref{hypothesesarefulfilledone} and explicitly check that $\TMF_1(n)_p$ is not real using \Cref{ktheorycalculation}. Also, the $2$-morphisms of \Cref{exampleofinversefixed} do not work here as $[-1]$ acts nontrivially on the point $p$ of order $n$, so long as $n\geq 3$. The same goes for the $2$-fold cover of $\M_{\Tate,\Z_p}$ and $\M_{\G_m,\Z_p}$, so $\Spf \Z_p\llbracket q\rrbracket$ and $\Spf \Z_p$, respectively, as their associated $\E_\infty$-rings are $\KU\llbracket q\rrbracket_p$ and $\KU_p$ which are not real.
\end{nonexample}

Let us clarify the origins of these natural families of $\E_\infty$-rings above.

\begin{mydef}
Fix an odd prime $p$. By \cite[\textsection12]{tmfbook} or \cite{hilllawson}, there is a sheaf $\O^\top$ of $\E_\infty$-rings on the small \'{e}tale site of $\M_{\Ell,\Z_p}$ which is uniquely defined (up to $1$-homotopy) by the fact that its affine sections come equipped with the natural structure of an elliptic cohomology theory; see \cite{uniqueotop}. In particular, given a formal Deligne--Mumford stack $\x$ with an \'{e}tale map $f\colon \x\to {\M}_{\Ell,\Z_p}$ there is an associated $\E_\infty$-ring $\O^\top(\x)$. Furthermore, as discussed in \cite[\textsection A]{tylerniko} and \cite[\textsection5]{hilllawson}, respectively (or in \cite[Prs.1.11 \& 1.18]{adamsontmf}), there are similarly defined sheaves of $\E_\infty$-rings $\O^\mult$ and $\O^\Tate$ on the moduli stacks $\M_{\G_m,\Z_p}$ and $\M_{\Tate,\Z_p}$.
\end{mydef}

Examples of sections of the above sheaves are the familiar (and also real, by \Cref{hypothesesarefulfilledone}) $\E_\infty$-rings
\[\O^\top(\M_{\Ell,\Z_p})=\Tmf_p\qquad \O^\top(\M^\sm_0(n)_{\Z_p})=\TMF_0(n)_p\]
\[\O^\Tate(\M_{\Tate,\Z_p})=\KO\llbracket q\rrbracket_p\qquad \O^\Tate(\M_{\Tate,\Z_p}^\sm)=\smKO_p\qquad \O^\mult(\M_{\G_m,\Z_p})=\KO_p\]
where $\M^\sm_0(n)$ denotes the moduli of smooth elliptic curves equipped with a cyclic subgroup of order $n$.

\begin{theorem}\label{hypothesesarefulfilledone}
Suppose $\x$ is a formal Deligne--Mumford stack and that we are given an \textbf{affine} morphism $f\colon \x\to \M_{\Ell,\Z_p}$. Make one of the following four additional assumptions:
\begin{enumerate}
\item Suppose that $\x=\Spf R$ is affine and $f$ is flat, leading us to a Landweber exact homotopy commutative elliptic cohomology theory $E$ associated with $f$.
\item Suppose that $f$ is \'{e}tale, which yields an $\E_\infty$-ring $E=\O^\top(\x)$.
\item Suppose that $f$ admits an \'{e}tale factorisation through $\M_\Tate$, which yields an $\E_\infty$-ring $\O^\Tate(\x)$.
\item Suppose that $f$ admits an \'{e}tale factorisation through $\M_{\G_m}$, which yields an $\E_\infty$-ring $\O^\mult(\x)$.
\end{enumerate}
If $f$ is inverse fixed, then $E$ is real with respect to the group $\{\pm 1\}$.
\end{theorem}

The above theorem is key for us to apply \Cref{ghtheorymain}, which in turn comes down to a $K$-theory calculation. See \cite[\textsection12.1 \& 12.5]{tmfbook} for the following definitions.

\begin{mydef}\label{definitionsoflevelstructures}
Write $\M^\ord_\Ell$ for the \emph{moduli stack of generalised elliptic curves with ordinary reduction over $p$-complete rings} and $\M^\ord_\Ell(p^\infty)$ for the continuous $\Z_p^\times$-torsor over $\M^\ord_\Ell$ defined by pairs $(C,\al)$ where $C$ is a generalised elliptic curve with ordinary reduction over a $p$-complete ring $R$ and $\al\colon \widehat{C}\simeq\widehat{\G}_m$ is an isomorphism of formal groups over $R$. From the lemma following the proof of \cite[Lm.12.5.1]{tmfbook}, the stack $\M^\ord_\Ell(p^\infty)$ is affine with global sections \emph{Katz' ring of $p$-adic modular forms $V$}.
\end{mydef}

\begin{prop}\label{ktheorycalculation}
Let $f\colon \x\to {\M}_{\Ell,\Z_p}$ and $E$ be as in hypotheses 1-4 of \Cref{hypothesesarefulfilledone} (here we do \textbf{not} assume $f$ to be inverse fixed). Let us write $f^\ord\colon \x^\ord\to \M^\ord_\Ell$ for the base-change of $f$ over $\M^\ord_\Ell\to {\M}_{\Ell,\Z_p}$. Then the zeroth $\K(1)$-local $\KU$-theory of $E$ fits into the following Cartesian diagram of formal stacks:
\begin{equation}\label{ktheoreycalculationsquare}\begin{tikzcd}
{\Spf \zK E}\ar[d]\ar[r]	&   {\M_\Ell^\ord(p^\infty)}\ar[d]  \\
{\x^\ord}\ar[r, "{f^\ord}"]			&   {\M^\ord_\Ell}
\end{tikzcd}\end{equation}
Moreover, the natural map of graded $\theta$-algebras $\K_\ast\otimes \zK E\to \cK E$ is an isomorphism.
\end{prop}


The above is a generalisation of \cite[Pr.12.6.1]{tmfbook} from affine schemes to general $\x$.

\begin{proof}
Under hypothesis $1$, this is precisely \cite[Pr.12.6.1]{tmfbook}. Continuing now under hypotheses 2-4, write $\x^\ord(p)$ be the formal Deligne--Mumford stack defined by the Cartesian diagram
\[\begin{tikzcd}
{\x^\ord(p)}\ar[r]\ar[d]		&	{\M_\Ell^\ord(p)}\ar[d]	\\
{\x^\ord}\ar[r, "{f^\ord}"]	&	{\M_\Ell^\ord}
\end{tikzcd}\]
where $\M_\Ell^\ord(p)$ is the $\F_p^\times$-\'{e}tale torsor (see \cite[Lm.12.5.1]{tmfbook}) of $\M_\Ell^\ord$ representing ordinary generalised elliptic curves curves $C$ with a choice of isomorphism $\mu_p\simeq \widehat{C}[p]$ of finite group schemes. By \cite[Lm.12.5.2]{tmfbook}, the stack $\M_\Ell^\ord(p)$ is affine, and $f^\ord$ is affine by base-change, hence $\x^\ord(p)$ is an affine formal Deligne--Mumford stack $\Spf W$. By evaluating either $\O^\top_{\K(1)}$ or $\O^\Tate_{\K(1)}$ on $\x^\ord(p)$, we obtain an $\E_\infty$-ring $E(p)$ and by descent for these \'{e}tale sheaves, an $\F_p^\times$-Galois extension of $\E_\infty$-rings $E\to E(p)$. In particular, this map induces an equivalence $E\simeq E(p)^{h\F_p^\times}$. By \cite[Pr.12.6.1]{tmfbook}, this proposition holds for $E(p)$, so in particular we now have the commutative diagram of formal Deligne--Mumford stacks
\[\begin{tikzcd}
{\Spf \zK E(p)}\ar[r, "{\pi'}"]\ar[d, "{\ga(p)}"]	&	{\Spf \zK E}\ar[d, "{\ga}"]\ar[r]	&	{\M_\Ell^\ord(p^\infty)}\ar[d]\\
{\Spf W=\x^\ord(p)}\ar[r, "{\pi}"]			&	{\x^\ord}\ar[r]				&	{\M_\Ell^\ord}
\end{tikzcd}\]
where we now know the outer rectangle is Cartesian. To see the right square is Cartesian and finish the proof, it suffices to show that the map $\ga$ is the $\F_p^\times$ quotient of the map $\ga(p)$. We know that $\pi$ is an $\F_p^\times$-quotient as it was constructed as the base-change of an $\F_p^\times$-torsor, so by naturality of the vertical maps, it suffices to show that $\pi'$ is also an $\F_p^\times$-quotient, in other words, that $\zK E\to \zK E(p)$ is the inclusion of $\F_p^\times$-fixed points. To see this consider the following chain of natural isomorphisms:
\[\cK E\xrightarrow{\simeq} \pi_\ast(\K\otimes E(p)^{h\F_p^\times}) \xleftarrow{\simeq} \pi_\ast((\K\otimes E(p))^{h\F_p^\times} \xrightarrow{\simeq} (\cK E(p))^{\F_p^\times}\]
The first follows as $E\to E(p)$ is an $\F_p^\times$-Galois extension, the second as $\K(1)$-local ambidexterity naturally identifies homotopy fixed points with homotopy orbits and hence homotopy fixed points commute with tensoring in one variable, and the third as the $E_2$-page of the $\F_p^\times$-HFPSS for $\K\otimes E(p)$ is concentrated in the zeroth row as $|\F_p^\times|=p-1$ is invertible. \\

For the ``moreover'' statement, the above shows that vertical morphisms in the diagram of $\theta$-algebras
\[\begin{tikzcd}
{\K_\ast \otimes \zK E}\ar[r]\ar[d]	&	{\cK E}\ar[d]	\\
{\K_\ast \otimes \zK E(p)}\ar[r]		&	{\cK E(p)}
\end{tikzcd}\]
are inclusions of $\F_p^\times$-fixed points. Our desired conclusion follows from the fact that the lower-horizontal morphism above is an isomorphism, a consequence of the fact that $E(p)$ is Landweber exact.
\end{proof}

\begin{proof}[Proof of \Cref{hypothesesarefulfilledone}]
To see that $E$ is real, notice that \Cref{ktheorycalculation} shows that $\cK E$ is concentrated in even degrees. Moreover, the facts that $\K_\ast$ is torsion-free and that $\zK E$ is flat over $\Spf \Z_p$ (and hence also flat over $\Spec \Z$) as all of the maps in (\ref{ktheoreycalculationsquare}) are flat and $\M^\ord_\Ell\to \Spf \Z_p$ is flat (as $\Spf V\to \M^\ord_\Ell$ is a $\Z_p^\times$-torsor and $V$ is flat over $\Z_p$) show that $\cK E$ is torsion-free. We are left to show that $\zK E$ has trivial $\{\pm1\}$-action---this part of the proof is inspired by the proof of \cite[Lm.12.7.14]{tmfbook} and \cite[Ex.6.12]{tmfwls}. \Cref{ktheorycalculation} shows that $\Spf \zK E$ represents pairs $(h,\al)$ where $h$ is a map of stacks $h\colon T\to \x^\ord$ from our test stack $T$ and $\al$ is an isomorphism of formal groups $\al\colon \widehat{\G}_m \simeq \widehat{C}_h$ over $T$, where $C_h$ is the elliptic curve defined by the composition $T\to \x^\ord\to \M_\Ell^\ord$. The $\Z_p^\times$-action on these $T$-points are given by sending the pair $(h,\al)$ to $(h,\al\circ [k])$ for each $p$-adic unit $k$, where $[k]\colon \widehat{\G}_m\to \widehat{\G}_m$ is the $k$-fold multiplication map on $\widehat{\G}_m$. This defines our $\Z_p^\times$-action on $\Spf \zK E$ in the slice $2$-category of stacks over $\M_\Ell$. To finish the proof, it suffices to show that the restricted $\{\pm 1\}$-action on $\zK E$ is trivial.\\

Notice that this restricted $\{\pm 1\}$-action on $\zK E$ is isomorphic in the slice $2$-category of stacks over $\M_\Ell$ to the $\{\pm 1\}$ action given by the identity on $\Spf \zK E$ as the inversion map $[-1]$ on the generalised elliptic curve $C_K$ defined over $\Spf \zK E$---this follows as the formal completion of this inversion map on $C_K$ induces the inversion map on the associated formal group $\widehat{C}_K$, and the isomorphism $\al$ is a homomorphism of formal groups and hence commutes with this inverse map. However, we assumed that $f$, and hence also $\x^\ord\to \x\to \M_{\Ell,\Z_p}$, is inverse fixed, meaning this $\{\pm1\}$-action is trivial. This finishes our proof.
\end{proof}

Although \Cref{hypothesesarefulfilledone} allows us to abstractly apply \Cref{ghtheorymain} to $\E_\infty$-rings coming from $\O^\top$, $\O^\Tate$, or $\O^\mult$, to obtain workable calculations, ie, to show certain vanishing results, we will need to know the following statement about the $\KO$-homology of these $\E_\infty$-rings.

\begin{prop}\label{hypothesesarefulfilledtwo}
Let $f$, $\x$, and $E$ be as in parts 2-4 of \Cref{hypothesesarefulfilledone}, assume $f$ is inverse fixed, and write $\Ga=\Z_p^\times/\{\pm 1\}$. Then in degree zero $\cKO E=A_\ast$ is ind-\'{e}tale over its $\Ga$-fixed points $A_0^0$, ${A}_0^0$ is the $p$-adic completion of a formally smooth algebra over $\Z_p$, and the continuous $\Ga$-cohomology groups $H^s_c(\Ga,A_\ast)$ vanish for $s\geq 1$.
\end{prop}

Our proof will essentially boil down to complex $K$-theory calculations and known facts.

\begin{proof}
For the first two statements, note that as $E$ is real with respect to $\{\pm1\}$ by \Cref{hypothesesarefulfilledone}, the complexification map
\[A_0=\zKO E\xrightarrow{\simeq} \zK E=B\]
is an isomorphism by \Cref{bottperiodicspectragivegood}. This reduces us to show that $B$ is ind-\'{e}tale over its $H=\Z_p^\times$-fixed points $B^H$, and that $B^H$ is the $p$-completion of a formally smooth algebra over $\Z_p$. As in the proof of \Cref{ktheorycalculation}, let us define stacks $\x^\ord(p)$ as $\F_p^\times$-Galois extensions of $\x^\ord$ with associated $\E_\infty$-ring $E(p)$. By \Cref{ktheorycalculation}, we also obtain the zeroth complex $K$-theory of $E(p)$ as an $\F_p^\times$-Galois extension of $B$, which we will denote as $B(p)$. Another application of \Cref{ktheorycalculation}, this time applied to $\x$, and $\y$, where  $\y$ is either $\M_\Ell$, $\M_\Tate$, or $\M_{\G_m}$, depending on which part of \Cref{hypothesesarefulfilledone} we find ourselves in, leads to the commutative diagram of formal Deligne--Mumford stacks
\begin{equation}\label{bigdiagramwithmycovertoo}\begin{tikzcd}
{\Spf B(p)}\ar[r]\ar[d]	&    {\Spf B}\ar[d]\ar[r, "{f^\ord(p^\infty)}"]	        &    {\y^\ord(p^\infty)}\ar[r]\ar[d]      &   {\M_\Ell^\ord(p^\infty)=\Spf V}\ar[d]  \\
{\x^\ord(p)}\ar[r]	&    {\x^\ord}\ar[r, "{f^\ord}"]	    &   {\y^\ord}\ar[r, "{g}"]      &   {\M^\ord_\Ell}
    \end{tikzcd}\end{equation}
in which all squares are Cartesian; $\y^\ord(p^\infty)$ is defined such that the right square is Cartesian. In particular, $\y^\ord(p^\infty)=\Spf W$ is affine as $g$ is affine. From this, we obtain the commutative diagram of $p$-adic rings
\[\begin{tikzcd}
{\Z_p}\ar[r]	&	{W^H}\ar[r]\ar[d] 			&   {B^H}\ar[d]\ar[r]	&	{B(p)^H}\ar[d] \\
		&	{W}\ar[r, "{f^\ord(p^\infty)}"]  	&   {B}\ar[r]		&	{B(p)}
\end{tikzcd}\]
where the vertical morphisms are inclusions of $H$-fixed points and the right horizontal morphisms are inclusions of $\F_p^\times$-fixed points. As $\x^\ord(p)$ is affine by construction, and by \cite[Cor.5.8(2)]{lawsonnaumannbp} we have $\x^\ord(p)\simeq \Spf B(p)^H$, then it follows by (\ref{bigdiagramwithmycovertoo}) and base-change that $B(p)^H\to B$ is ind-\'{e}tale. By taking $\F_p^\times$-fixed points, it follows that $B^H\to B$ is also an ind-\'{e}tale extension. Similarly, to see $B^H$ is the $p$-completion of a formally smooth $\Z_p$-algebra, we use the fact that $B\to B(p)$ is an \'{e}tale cover and the fact that formally smooth is an adjective which is \'{e}tale local on the source; see \cite[\href{https://stacks.math.columbia.edu/tag/061K}{061K}]{stacks}. To see that $B(p)^H$ is a $p$-completion of a formally smooth $\Z_p$-algebra, we appeal to \cite[Lm.5.5(3)]{lawsonnaumannbp}.\\

For the cohomological vanishing, first note the isomorphism of $\Ga$-modules
\[A_{4t}=\ftKO E\simeq \zKO E(\chi_{2t})=A_0(\chi_{2t})\]
where $\chi_{t}\colon \Ga\to \Ga$ is the continuous character sending $k\mapsto k^{t}$. This comes from the action of the stable Adams operations on $\KO_\ast$ and the fact that $E$ is real by \Cref{hypothesesarefulfilledone}. We now use a general fact that if $R\to S$ is an ind-Galois extension with respect to $\Ga$, then the continuous cohomology groups $H^s_c(\Ga, S(\chi))$ vanish for all characters $\chi\colon \Ga\to \Ga$; this is explained clearly in the proof of part 2 of \cite[Lm.5.5]{lawsonnaumannbp}. This finishes the proof, as we know that $A_0^0\to A_0$ is a $\Ga$-Galois extension, hence so is $A_\ast^0\to A_\ast$.
\end{proof}

Both \Cref{hypothesesarefulfilledone} and \Cref{hypothesesarefulfilledtwo} have refinements to $G$-homotopy fixed points of various $\E_\infty$-rings.

\begin{cor}\label{corollarythatallisfufilled}
Let $G$ be a nontrivial finite subgroup of $\Z_p^\times$. Then $\Tmf_p^{hG}$, $\TMF_p^{hG}$, $\KO\llbracket q\rrbracket_p^{hG}$, and $\smKO_p^{hG}$ are all real, where all $G$-actions are induced by the $p$-adic stable Adams operations. Moreover, if $E$ is any of these $\E_\infty$-rings, then in degree zero $\cL E=A_\ast$ is ind-\'{e}tale over its $\Ga$-fixed points $A_0^0$, ${A}_0^0$ is the $p$-adic completion of a formally smooth algebra over $\Z_p$, and the continuous $\Ga$-cohomology groups $H^s_c(\Ga,A_\ast)$ vanish for $s\geq 1$.
\end{cor}

\begin{proof}
The classical calculations that $\zK \Tmf\simeq V$ and $\zK \K\llbracket q\rrbracket\simeq \Cont(\Z_p^\times,\Z_p)$, where $V$ is the $p$-completion of Katz' ring of divided congruences (\cite[\textsection 1.4]{handbooktmf}) yield the following equivalences
\[\zK\Tmf^{hG}\simeq (\zK \Tmf)^{G}\simeq V^G\qquad \zK \K\llbracket q\rrbracket^G\simeq \Cont(\Z_p^\times,\Z_p)^G=\Cont(\Ga,\Z_p)\]
as the $p$-adic Adams operations on these $\E_\infty$-rings, so \cite[Th.A]{adamsontmf} for $\Tmf_p$, \cite[\textsection5.5]{luriestheorem} for the rest, induce the $\psi$-module $\Z_p^\times$-action on the above groups. It is clear from the above descriptions that $G$ acts trivially on the groups above, and the smooth cases, so $\TMF_p$ and $\smK_p$, are similar. This combined with the computations
\[\cK E=\cK R^{hG}\simeq (\cK R)^G\simeq \K_\ast\otimes (\zK R)^G\simeq \K_\ast\otimes \zK E\]
for any of the $E=R^{hG}$ above, shows that these $\E_\infty$-rings are real with respect to $G$. For the ``moreover'' statement, one simply copies the proof for \Cref{hypothesesarefulfilledtwo}, only needing to replace $\{\pm 1\}$ with $G$ and $\KO$ with $\L$.
\end{proof}


\subsection{Uniqueness of $p$-adic topological $q$-expansion map (\Cref{mainuniquenessstatement,mainuniquenessstatementpgrowth})}\label{uniquenessofqexpasniosection}

In this section, we prove \Cref{mainuniquenessstatement,mainuniquenessstatementpgrowth}. First, we have the following $p$-adic generalisation of \Cref{mainuniquenessstatement}.

\begin{theorem}\label{generalisedmainuniquenessstatmet}
Let $p$ be an odd prime, $f,E$ be as in parts 2-4 of \Cref{hypothesesarefulfilledone} and $f',E'$ be as in parts 3-4 of \Cref{hypothesesarefulfilledone}. Suppose that both $f$ and $f'$ are inverse fixed. Then every map of $\E_\infty$-rings $f\colon E\to E'$ is uniquely determined by its effect on the zeroth $\K(1)$-local $\KO$-homology (or $\KU$-homology) group up to $3$-homotopy.
\end{theorem}

To prove the above theorem, we will use \Cref{ghtheorymain} together with the following vanishing statement of Andr\'{e}--Quillen cohomology groups.

\begin{prop}\label{vanishingofobstructions}
Fix an odd prime $p$ and a nontrivial finite subgroup $G\leq \Z_p^\times$. Suppose $A_\ast$ is a graded $\bar{\theta}$-algebra and that $M_\ast$ is a graded $\bar{\theta}$-module. Write $A_\ast^0=A_\ast^\Ga$ where $\Ga=\Z_p^\times/G$. Suppose the following conditions hold:
\begin{enumerate}
    \item Both $A_\ast$ and $M_\ast$ are both torsion-free and real with respect to $G$.
    \item $A_0^0$ is formally smooth over $\Z_p$.
    \item The continuous group cohomology groups $H_c^s(\Ga,M_u)$ vanish for $s>0$ and every $u\in \Z$.
    \item $A_0$ is ind-\'{e}tale over $A_0^0$.
\end{enumerate}
Then $H^s_{\bar{\theta}}(A_\ast/\L_\ast, M_\ast[t])=0$ for $s\geq 2$ and all $t$ not divisible by $2g$.
\end{prop}

To prove this vanishing statement, one would simply copy Lawson--Naumann's proof of \cite[Lm.5.15]{lawsonnaumannbp}, replacing their $\Z_p^\times$ with our $\Ga$; also see Behrens' proof of \cite[Lm.12.7.5 \& Lm.12.7.13]{tmfbook} for a similar argument. Let us avoid too much repetition and forgo proof. The vanishing for $t$ not divisible by $2g$ comes from the fact that $M_\ast$ is real, so $L_\ast\otimes M_0=M_\ast$, and $L_\ast$ is concentrated in degrees divisible by $2g$.

\begin{proof}[Proof of \Cref{generalisedmainuniquenessstatmet}]
First note that the zeroth $\KO$- and $\KU$-homology of real spectra with respect to $\{\pm1\}$ are naturally isomorphic by \Cref{bottperiodicspectragivegood}.\\

Now, consider the classical Goerss--Hopkins obstruction theory of \cite[Cor.4.4]{gh04} based on $\KU$. By \Cref{hypothesesarefulfilledone}, we see the spectra in sight are real and \Cref{cohomologyiseasyforrealguys} states it suffices to compute Andr\'{e}--Quillen cohomology of their $\KO$-homology. By \Cref{hypothesesarefulfilledtwo} and \Cref{vanishingofobstructions} (for $G=\{\pm 1\}$) we see that enough obstruction groups vanish to conclude that $f$ is uniquely determined up to $1$-homotopy---this fact is already well-known.\\

By \Cref{reductiontohomotopygroups}, it now suffices to show that the groups $\pi_k(\CAlg(E,E'),f)$ vanish for all $0 < k<3$. Note that both $E$ and $E'$ are $p$-complete and $E'$ is $\K(1)$-local by assumption, and that according to \Cref{hypothesesarefulfilledone} both $E$ and $E'$ are real. We may then use \Cref{ghtheorymain} to obtain a spectral sequence that abuts to the desired homotopy groups:
\[H^s_{\bar{\theta}}(\cKO E/\KO_\ast, \cKO E'[-t])\Longrightarrow \pi_{t-s}(\CAlg(E,E'),f)\]
Combining \Cref{hypothesesarefulfilledtwo} with \Cref{vanishingofobstructions}, we see these groups vanish for $s\geq 2$ and $t$ not divisible by $4$. This spectral sequence above is concentrated in the $s=0$ and $s=1$ rows in degrees $t$ divisible by $4$. This immediately yields the desired vanishing of $\pi_i$ for $0<i<3$.
\end{proof}

\begin{remark}\label{failiureoftwopimraryremark}
At the prime $2$, all of the arguments above work as in the odd primary case, except that the vanishing result \Cref{vanishingofobstructions} at the prime $2$ (see \cite[Lm.12.7.13]{tmfbook}) does not guarantee the desired vanishing of obstruction groups. In essence, this comes down to the fact that $\pi_1\KO\neq 0$.
\end{remark}

To prove \Cref{mainuniquenessstatement}, we combine the general $p$-adic uniqueness statement \Cref{generalisedmainuniquenessstatmet} with standard rational stable homotopy theory.

\begin{proof}[Proof of \Cref{mainuniquenessstatement}]
Consider the Cartesian diagram of spaces
\[\begin{tikzcd}
{\CAlg(\tmf,\KO\llbracket q\rrbracket[\frac{1}{2}])}\ar[r]\ar[d]	&	{\CAlg(\tmf, \prod \KO\llbracket q\rrbracket_p)}\ar[d]	\\
{\CAlg(\tmf,\KO\llbracket q\rrbracket_\Q)}\ar[r]			&	{\CAlg\left(\tmf, \left(\prod \KO\llbracket q\rrbracket_p\right)_\Q\right)}
\end{tikzcd}\]
induced by the arithmetic fracture square for $\KO\llbracket q\rrbracket[\frac{1}{2}]$, where the products are taken over all odd primes. Applying the functor $\tau_{\leq 2}$ does not generally commute with limits of spaces, however, from the long exact sequence of a fibre product on homotopy groups, we see that the diagram
\begin{equation}\label{littlecartesianfiend}\begin{tikzcd}
{\tau_{\leq 2}\CAlg(\tmf,\KO\llbracket q\rrbracket[\frac{1}{2}])}\ar[r]\ar[d]	&	{\tau_{\leq 2}\CAlg(\tmf, \prod_{p\neq 2}\KO\llbracket q\rrbracket_p)}\ar[d]	\\
{\tau_{\leq 2}\CAlg(\tmf,\KO\llbracket q\rrbracket_\Q)}\ar[r]			&	{\tau_{\leq 2}\CAlg\left(\tmf, \left(\prod_{p\neq 2}\KO\llbracket q\rrbracket_p\right)_\Q\right)}
\end{tikzcd}\end{equation}
is Cartesian if and only if the induced map
\[\pi_3\CAlg\left(\tmf, \left(\prod_{p\neq 2}\KO\llbracket q\rrbracket_p\right)_\Q\right)\to \pi_2 \CAlg(\tmf,\KO\llbracket q\rrbracket[\frac{1}{2}])\]
is zero. We claim that the $\pi_3$ above vanishes, which would imply that (\ref{littlecartesianfiend}) is Cartesian. To see the above $\pi_3$ vanishes, note the classical fact that $\tmf_\Q$ is the free formal cdga on two variables $x$ and $y$, where $|x|=8$ and $|y|=12$; see \cite[Pr.4.47]{hilllawson}. This implies that for a rational $\E_\infty$-ring $R$, the mapping space $\CAlg(\tmf,R)$ is naturally equivalent to $\Omega^{8}R\times \Omega^{12}R$. In particular, as the homotopy groups of $(\prod \KO\llbracket q\rrbracket_p)_\Q$ are concentrated in degrees divisible by $4$, we see that the above $\pi_3$ vanishes.\\

We claim that all of the spaces in (\ref{littlecartesianfiend}) are discrete. If this were true, then a morphism of $\E_\infty$-rings from $\tmf$ to $\KO\llbracket q\rrbracket[\frac{1}{2}]$ is uniquely determined by its combined effect upon taking $p$-completions for all odd primes $p$ and rationalisation up to $3$-homotopy. Using the above-mentioned fact that $\tmf_\Q$ is a free formal rational cdga and the fact that the homotopy groups of $\KO\llbracket q\rrbracket_\Q$ and $(\prod \KO\llbracket q\rrbracket_p)_\Q$ are supported in degrees divisible by $4$, we see that the lower two spaces in (\ref{littlecartesianfiend}) are discrete. Moreover, we note that the effect of a morphism of $\E_\infty$-rings from $\tmf$ to $\KO\llbracket q\rrbracket [\frac{1}{2}]$ on rational homotopy groups is determined by its effect on zeroth rational $\KU$-homology. Indeed, this is because rationally we have an equivalent of spectra $\KU_\Q\simeq \bigoplus \Sph_\Q[2n]$ and all above rational $\E_\infty$-rings $R$ have homotopy groups supported in even degrees, hence there is an identification $\pi_\ast R\simeq \K_0 R$, natural in $R$.\\

Noting that $\tmf\to \Tmf$ is a $\K(1)$-local equivalence (\cite[p.30]{handbooktmf}), we see that the upper-right space of (\ref{littlecartesianfiend}) is discrete by \Cref{reductiontohomotopygroups} and \Cref{generalisedmainuniquenessstatmet} with $\x=\M_{\Ell,\Z_p}$ and $\y=\M_{\Tate,\Z_p}$.\\

In summary, we see that all of the spaces in (\ref{littlecartesianfiend}) are discrete, and that a morphism of $\E_\infty$-rings $\tmf\to \KO\llbracket q\rrbracket[\frac{1}{2}]$ is uniquely determined by its effect on zeroth $p$-adic $\KU$-homology for all odd primes $p$ and zeroth rational $\KU$-homology up to $3$-homotopy. As all of these homology groups are completions or localisations of the zeroth integral $\KU$-homology group, we see such morphisms of $\E_\infty$-rings are uniquely determined by their effect on zeroth $\KU$-homology up to $3$-homotopy.
\end{proof}

To prove \Cref{mainuniquenessstatementpgrowth}, we follow the general outline of the proof of \Cref{generalisedmainuniquenessstatmet}.

\begin{proof}[Proof of \Cref{mainuniquenessstatementpgrowth}]
As the zeroth $\L$- and $\K$-homology of real spectra are naturally isomorphic by \Cref{bottperiodicspectragivegood}, then the classical Goerss--Hopkins obstruction theory of \cite[Cor.4.4]{gh04} based on $\KU$, combined with \Cref{cohomologyiseasyforrealguys} and \Cref{corollarythatallisfufilled} show we are left to compute Andr\'{e}--Quillen cohomology of $\L$-homologies as $\bar{\theta}$-algebras. By \Cref{corollarythatallisfufilled} and \Cref{vanishingofobstructions}, enough obstruction groups vanish to show that $f$ is uniquely determined up to $1$-homotopy.\\

By \Cref{reductiontohomotopygroups}, it now suffices to show that the groups $\pi_k(\CAlg(\tmf_p^{hG},\K\llbracket q\rrbracket_p^{hG}),f)$ vanish for all $0 < k<2g-1$. By \Cref{corollarythatallisfufilled}, we see both $\E_\infty$-rings are real, so we can apply \Cref{ghtheorymain} to obtain a spectral sequence that abuts to the desired homotopy groups:
\[H^s_{\bar{\theta}}(\cL \tmf^{hG}/\L_\ast, \cL \K\llbracket q\rrbracket^{hG}[-t])\Longrightarrow \pi_{t-s}(\CAlg(\tmf^{hG},\K\llbracket q\rrbracket^{hG}),f)\]
Combining \Cref{corollarythatallisfufilled} with \Cref{vanishingofobstructions}, we see these groups vanish for $s\geq 2$ and $t$ not divisible by $2g$, and we are done.
\end{proof}


\subsection{Higher functoriality of elliptic Adams operations (\Cref{adamsoperationshomotopiesintegral,adamsoperationshomotopiespadic})}\label{adamsoperationssection}

In \cite{adamsontmf}, we defined morphisms of $\E_\infty$-rings $\psi^k\colon \tmf[\frac{1}{k}]\to \tmf[\frac{1}{k}]$ for each integer $k$ which we called \emph{stable Adams operations} due to either their construction, effect on homotopy groups, or relationship with the classical operations on topological $K$-theory. Our construction of these operations used Goerss--Hopkins obstruction theory and an unfortunate consequence of these techniques was a failure to prove any relationship between these operations for varying $k$. Using the homotopical uniqueness of the $q$-expansion map from the previous section, we will see that these Adams operations on $\Tmf$ satisfy the usual relations $\psi^{k\ell}\simeq \psi^k\psi^\ell$ up to $2$-homotopy.\\

Let us start with a proof of \Cref{adamsoperationshomotopiespadic}---actually, we will first prove a slight generalisation for dualisable topological modular forms.

\begin{theorem}\label{adamsoperationshomotopies}
Let $p$ be an odd prime and $G\leq \Z_p^\times$ be a nontrivial finite subgroup with order $g$. Writing $\Ga=\Z_p^\times/G$, there exists a functor
\[\Psi_{p,G}\colon B\Ga\to \h_{2g-2}\CAlg\]
sending the unique point of $B\Ga$ to $\Tmf_p^{hG}$ and each $p$-adic unit $\overline{k}$ to $\psi^k$. 
\end{theorem}

The argument above is \textbf{not} that the Adams operations $\psi^k$ on $\Tmf_p$ exist and are highly unique by obstruction theory, but rather we take highly coherent homotopies between Adams operations on $\TMF_p$, using their modular descriptions, and then glue these homotopies together using the aforementioned obstruction theory.

\begin{proof}
To start with, let us construct a functor $\Psi_1\colon B\Ga\to \h_1\CAlg$ with the desired properties. This is equivalent to constructing a morphism of $\E_1$-monoids $\Ga\to \tau_{\leq 0}\CAlg(\Tmf_p,\Tmf_p)$. To do this, consider the Cartesian diagrams of $\E_\infty$-rings
\[\begin{tikzcd}
{\Tmf_p}\ar[r]\ar[d]	&	{\KO\llbracket q\rrbracket_p}\ar[d]	\\
{\TMF_p}\ar[r]		&	{\smKO_p}
\end{tikzcd}\qquad\qquad
\begin{tikzcd}
{\Tmf_p^{hG}}\ar[r]\ar[d]	&	{\KO\llbracket q\rrbracket_p^{hG}}\ar[d]	\\
{\TMF_p^{hG}}\ar[r]		&	{\smKO_p^{hG}}
\end{tikzcd}\]
given by taking global sections of \cite[Df.5.10]{hilllawson} and then taking $G$-homotopy fixed points, which induces the following diagram of spaces:
\begin{equation}\label{diagramoftruncations}\begin{tikzcd}
{\tau_{\leq d}\CAlg(\Tmf_p^{hG},\Tmf_p^{hG})}\ar[r]\ar[d]  	& {\tau_{\leq d}\CAlg(\Tmf_p^{hG},\KO\llbracket q\rrbracket_p^{hG})}\ar[d] \\
{\tau_{\leq d}\CAlg(\Tmf_p^{hG},\TMF_p^{hG})}\ar[r] 		& {\tau_{\leq d}\CAlg(\Tmf_p^{hG},\smKO_p^{hG})}
\end{tikzcd}\end{equation}
As in the proof of \Cref{mainuniquenessstatement} above, this above square would be Cartesian if we could show that $\pi_{d+1} \CAlg(\Tmf_p^{hG},\smKO_p^{h_G})$ vanishes. Let us now restrict to those $d$ in the range $0\leq d\leq 2g-3$. For $G=\{\pm 1\}$, we see that the desired $\pi_{d+1}$ vanishes by \Cref{generalisedmainuniquenessstatmet} (for $\x=\M_{\Ell,\Z_p}$ and $\y=\M_{\Tate,\Z_p}^\sm$) and \Cref{reductiontohomotopygroups}. For a general $G$, we appeal to the spectral sequence of \Cref{ghtheorymain} and the vanishing of much of the $E_2$-page from \Cref{corollarythatallisfufilled} and \Cref{vanishingofobstructions}, which together show this group vanishes. In particular, when $d=0$ the fact that (\ref{diagramoftruncations}) is Cartesian means that to construct a map of spaces $\Ga\to [\Tmf_p^{hG},\Tmf_p^{hG}]_{\CAlg}$, it suffices to construct two maps
\[\Psi_\T\colon\Ga\to [\Tmf_p^{hG},\TMF_p^{hG}]_{\CAlg}\qquad\qquad \Psi_{\K}\colon\Ga\to [\Tmf_p^{hG},\KO\llbracket q\rrbracket_p^{hG}]_{\CAlg}\]
which agree when restricted to $[\Tmf_p^{hG},\smKO_p^{hG}]_{\CAlg}$---we will come back to the $\E_1$-structures shortly. Define the first map $\Psi_\T$ as the composition
\[\Z_p^\times/\{\pm1\}\xrightarrow{\Psi^\sm_p} \CAlg(\TMF_p^{hG},\TMF_p^{hG})\to \CAlg(\Tmf_p^{hG},\TMF_p^{hG})\to \tau_{\leq 0}\CAlg(\Tmf_p^{hG},\TMF_p^{hG})\]
where $\Psi^\sm_p$ is a morphism of $\E_1$-monoids, a $p$-adic version of \cite[Pr.2.37]{heckeontmf}, itself a corollary of the fact that $\TMF_p$ is the global sections of the sheaf $\O^\top_\BTtwo$ on the moduli stack of elliptic curves; see \cite[Th.5.17]{luriestheorem}. The map $\Psi_{\K}$ can either be constructed similarly to $\Psi^\sm_p$ using a moduli construction, as in \cite[Pr.1.18]{adamsontmf}, or one can use the uniqueness statement \Cref{generalisedmainuniquenessstatmet} to see that morphisms from $\Tmf_p^{hG}$ to $\KO\llbracket q\rrbracket_p^{hG}$ are uniquely determined by their effect on $p$-adic $\K$-homology up to $(2g-1)$-homotopy, hence there exists a homotopy between $\psi^k\psi^\ell$ and $\psi^{k\ell}$, as these morphisms have the same modular description on $p$-adic $\K$-homology before taking $G$-homotopy fixed points; see \Cref{ktheorycalculation} for the $\K$-theory calculations and the proof of \Cref{hypothesesarefulfilledone} to see the modular description of these Adams operations. This second approach, using uniqueness up to homotopy and the modular description of Adams operations, shows that these two maps $\Psi_\T$ and $\Psi_{\K}$ agree when restricted to $\smKO_p^{hG}$. We then obtain a morphism of spaces
\begin{equation}\label{eonemonoidmap}\Psi_1\colon \Ga\to [\Tmf_p^{hG},\Tmf_p^{hG}]_{\CAlg}\end{equation}
and to check this is a morphism of discrete $\E_1$-monoids, so equivalently monoids of sets, we reuse the arguments from above: to check the relation between $\psi^k\psi^\ell$ and $\psi^{k\ell}$ holds inside the discrete monoid $[\Tmf_p^{hG},\Tmf_p^{hG}]_{\CAlg}$, we need to produce a homotopy between these morphisms. Using (\ref{diagramoftruncations}) now for $d=1$, we obtain the desired homotopy restricted to $\TMF_p^{hG}$ from $\Psi^\sm_p$ in the construction of $\Psi_\T$, and also when restricted to $\KO\llbracket q\rrbracket_p^{hG}$ using the homotopy uniqueness of morphisms from $\Tmf_p^{hG}$ to $\KO\llbracket q\rrbracket_p^{hG}$. These homotopies agree on $\smKO_p$ as such morphisms are unique up to $(2g-1)$-homotopy, hence we obtain a homotopy $\psi^k\psi^\ell\simeq \psi^{k\ell}$ on $\Tmf_p^{hG}$, and we see that (\ref{eonemonoidmap}) can be upgraded to a morphism of discrete $\E_1$-monoids.\\

To construct the functor promised in the theorem, let us first note that the natural map $\Tmf_p\to \TMF_p$ is a localisation at an element $\Delta^r\in \pi_{24r}\Tmf_p$, where for $p=3$ we have $r=3$, and else $r=1$. This persists to the $G$-homotopy fixed points. Indeed, it is clear from the effect on homotopy groups that the natural map $\Tmf_p^{hG}\to \TMF_p^{hG}$ is a localisation at $\Delta^s\in \pi_{24s}\Tmf_p^{hG}$, where $s$ is divisible by $2p-2$; smaller choices of $s$ are possible, but this does not play a role in our argument. This implies that the natural map of spaces (whose truncation appears in (\ref{diagramoftruncations}))
\[\CAlg(\Tmf_p^{hG},\Tmf_p^{hG})\to \CAlg(\Tmf_p^{hG},\TMF_p^{hG})\simeq \CAlg(\TMF_p^{hG},\TMF_p^{hG})\]
can be identified with the natural map of $\E_1$-monoids
\begin{equation}\label{firstmapofeonemonoids}\CAlg(\Tmf_p^{hG},\Tmf_p^{hG})\to \CAlg(\TMF_p^{hG},\TMF_p^{hG})\end{equation}
defined by inverting this $\Delta^{s}$. We now claim that for $0\leq d\leq 2g-3$, the following natural diagram of $\E_1$-monoids
\begin{equation}\label{finalcartdiagram}\begin{tikzcd}
{\tau_{\leq d}\CAlg(\Tmf_p^{hG},\Tmf_p^{hG})}\ar[r]\ar[d]	&	{\tau_{\leq d}\CAlg(\TMF_p^{hG},\TMF_p^{hG})}\ar[d]	\\
{\tau_{\leq 0}\CAlg(\Tmf_p^{hG},\Tmf_p^{hG})}\ar[r]		&	{\tau_{\leq 0}\CAlg(\TMF_p^{hG},\TMF_p^{hG})}\end{tikzcd}\end{equation}
is Cartesian. Indeed, this follows from the Cartesian diagram (\ref{diagramoftruncations}) and \Cref{mainuniquenessstatementpgrowth} which shows that the spaces in the right of (\ref{diagramoftruncations}) are discrete, hence the natural map (\ref{firstmapofeonemonoids}) is an isomorphism on homotopy groups in positive degrees. Therefore, to construct our desired map of $\E_1$-monoids $\Ga\to \tau_{\leq 2g-3}\CAlg(\Tmf_p^{hG},\Tmf_p^{hG})$ it suffices to construct maps of $\E_1$-monoids into the other factors of the Cartesian diagram (\ref{finalcartdiagram}). This is achieved using our map $\Psi_1$ and $\Psi^\sm_p$ discussed above, which agree on $\TMF_p^{hG}$ as the construction of $\Psi_1$ also used $\Psi^\sm_p$.
\end{proof}

\begin{proof}[Proof of \Cref{adamsoperationshomotopiespadic}]
Post-compose the functor of \Cref{adamsoperationshomotopies} with the connective cover functor.
\end{proof}

Finally, we can prove \Cref{adamsoperationshomotopiesintegral} by combining \Cref{adamsoperationshomotopiespadic} with standard rational information.

\begin{proof}[Proof of \Cref{adamsoperationshomotopiesintegral}]
Fix an integer $n$ with a collection of primes $p|n$ and write $M=\prod_{p|n}\N$. For each odd prime $\ell \nmid n$, the inclusion of monoids
\[M\xrightarrow{(p_i^{e_i})\mapsto \prod p_i^{e_i}} \N\to \Z\to \Z_\ell\]
factors as a injection $M\to \Z_p^\times$, and remains injective after a further quotient $M\to \Z_p^\times/\{\pm 1\}$ as $M$ contains no elements with additive inverses other than $0$. For all such odd $\ell$, we then have a morphism of $\E_1$-monoids
\[M\to\tau_{\leq 1}\CAlg(\tmf_\ell,\tmf_\ell)\]
by \Cref{adamsoperationshomotopiespadic} with respect to $G=\{\pm 1\}$. There is an analogous rational construction
\[\Psi_\Q^n\colon M\to \tau_{\leq 1}\CAlg(\tmf_\Q,\tmf_\Q)\]
which we define as follows: First, note that the rational Goerss--Hopkins obstruction theory in Step 1 of \cite[\textsection12.9]{tmfbook} shows that the codomain of $\Psi_\Q^n$ is discrete and equivalent to the set of maps of $\Q$-algebras endomorphisms of $\pi_\ast\tmf_\Q$. Next, recall from the proof of \Cref{mainuniquenessstatement} that $\tmf_\Q$ is the free formal cdga on two variables $x$ and $y$, where $|x|=8$ and $|y|=12$. We then define $\Psi_\Q^n$ by sending $p_i$ to the endomorphism of $\tmf_\Q$ sending $x$ to $p_i^4x$ and $y$ to $p_i^6y$ and extending multiplicatively.\\

To glue together $\Psi_\ell^n$ with $\Psi_\Q^n$, use the arithmetic fracture square for $\tmf[\frac{1}{2mn}]$ which yields the diagram of spaces
\[\begin{tikzcd}
{\tau_{\leq 1}\CAlg(\tmf[\frac{1}{2mn}],\tmf[\frac{1}{2mn}])}\ar[r]\ar[d]	&	{\prod_{2\neq \ell\nmid n}\tau_{\leq 1}\CAlg(\tmf_\ell,\tmf_\ell)}\ar[d]	\\
{\tau_{\leq 1}\CAlg(\tmf_\Q,\tmf_\Q)}\ar[r]				&	{\tau_{\leq 1}\CAlg(\tmf_\Q,(\prod_{2\neq \ell\nmid n}\tmf_\ell)_\Q)}
\end{tikzcd}\]
where the left vertical map and upper horizontal map are morphisms of $\E_1$-monoids. Again, as in the proof of \Cref{mainuniquenessstatement}, for any rational $\E_\infty$-ring $R$, the space $\CAlg(\tmf,R)$ is equivalent to $\Omega^{\infty+8}R\times\Omega^{\infty+12}R$, so the above diagram is Cartesian, as $\pi_2 \CAlg(\tmf_\Q,(\prod \tmf_\ell)_\Q)$ vanishes. The construction of $\Psi^n$ now comes down to showing that $\prod\Psi^n_\ell$ and $\Psi_\Q^n$ agree when restricted to the space in the lower-right corner. Again appealing to the rational Goerss--Hopkins obstruction theory in Step 1 of \cite[\textsection12.9]{tmfbook} shows that this lower-right space is discrete and equivalent to the associated mapping set of graded $\Q$-algebras after taking homotopy groups. In particular, the calculations of $\psi^k$ found in \cite[Th.C]{adamsontmf} show that both $\prod\Psi^n_\ell$ and $\Psi_\Q^n$ agree when restricted to this space in the lower-right corner, leading to our desired morphism of spaces
\[\Psi^n\colon M\to \tau_{\leq 1}\CAlg(\tmf[\frac{1}{2n}],\tmf[\frac{1}{2n}]).\]
Again, to see that this rectifies to a morphism of $\E_1$-monoids, we want to check that there exist homotopies recognising this fact. As in the proof of \Cref{adamsoperationshomotopies} above, one should do this by referring to the uniqueness of morphisms of $\E_\infty$-rings from $\tmf_\Q$ to $(\prod \tmf_\ell)_\Q$, and this is clear from the above characterisation of $\tmf_\Q$ and the fact that $(\prod \tmf_\ell)_\Q$ has trivial $\pi_i$ for $1\leq i\leq 7$.\\

Similarly, we can lift the homotopies $\psi^{-1}\simeq\id$ on $\tmf_\ell$ using \Cref{adamsoperationshomotopies} and on $\tmf_\Q$ (by construction), to a homotopy on $\tmf[\frac{1}{2}]$. This finishes the proof.
\end{proof}

\begin{remark}\label{fromtmftotmf}
The only reason we restricted from $\Tmf$ to $\tmf$ to prove \Cref{adamsoperationshomotopiesintegral} is as rationally $\tmf$ has a simple with universal property. Perhaps one can use the Cartesian square
\[\begin{tikzcd}
{\Tmf_\Q}\ar[r]\ar[d]	&	{\tmf_\Q[c_4^{-1}]}\ar[d]	\\
{\TMF_\Q}\ar[r]		&	{\TMF_\Q[c_4^{-1}]}
\end{tikzcd}\]
to prove a version of \Cref{adamsoperationshomotopiesintegral} for $\Tmf$.
\end{remark}


\addcontentsline{toc}{section}{References}


\bibliography{/Users/jackdavies/Dropbox/Work/references} 


\bibliographystyle{alpha}
\end{document}